# From Concrete to Abstract in Indian Mathematics

Jaidev Dasgupta

## ABSTRACT

Despite the extensive amount of scholarly work done on Indian mathematics in the last 200 years, the historical conditions under which it originated and evolved has not been studied much. The focus has been more on achievements than on how they developed. One also tends to read the ancient texts with the present concepts and methods in mind. The absence of writing over a long stretch of Indian history too gets overlooked in such readings. The purpose of this article is to explore the journey of mathematics by examining what the ancient texts on arithmetic, geometry and algebra tell us about the nature of mathematics in their times. These investigations reveal that over a period of a thousand or more years Indian mathematics transitioned from concrete and context-bound phase to context-free, abstract phase accompanied by several conceptual leaps. Invention of writing in the 3$^{rd}$ century BCE greatly facilitated this transition.

## 1. Introduction

It is well established now that ancient humans developed their sense of number through concrete, lived experience. Born with natural mathematical ability or "math instinct," a vague sense of numbers, there is a good likelihood that in different parts of the world in antiquity humans developed early stages of mathematics. Twenty to thirty thousand years old or even older archaeological finds such as bones and stones with scratch marks are indicative of record keeping of numerical information and, perhaps, of number counting in the remote past. Notched or tally sticks and human hands are also considered among the ancient tools for counting and calculations, though fingers could not be used for record keeping. Although very little is known about the origins of numbers and counting, a wide range of numeration systems have been invented by humanity since antiquity. The cognitive capacity of counting is believed to have developed in humans sometime between the Upper Paleolithic, pre-historic time and the historic time. The need to keep account of personal material possessions, e.g., jewelry, weapons, animals etc. or of the tally marks (scratches/notches) on bones and sticks is believed to have inspired the need for counting (Barras, 2021). Fast forward and we find ancient Egyptians counting oxen, goats and prisoners gained from the victory of the predynastic king Narmer over lower Egypt in the late fourth millennia BCE (Clagett, 1999; p. 2). Archaeological evidence from around the same time or even before show the use of clay tokens of multiple shapes and their assemblages as clay balls, bullae, by the Mesopotamians for recording goods such as wheat, cloth, oil, etc. and their quantities received or disbursed. These tokens or counters were concrete representations of specific items that believably led to the invention of cuneiform writing system representing objects and numbers (Schmandt-Besserat, 1996; Damerow, 2012).

On the Indian subcontinent, however, there is no archaeological evidence of tokens, buildings, or rock inscriptions that could demonstrate the knowledge of numbers and arithmetic in the period from the 18th century to 3rd century BCE. However, the ancient texts make extensive references to mathematics - especially in decimal enumeration, astronomy, calendrics, and sacrificial altar-making. The *Śulvasūtras* for altar-making and the astronomical text *Vedāṅga Jyotiṣa* for calendrics are prime examples of mathematical knowledge in antiquity. Much work has been done on the history of Indian mathematics in the last 200 years. But its focus is primarily on the past achievements, not the cultural conditions and society under which they were made. Exploring through ancient texts, this article makes a case for the concrete and context-dependent nature of Indian mathematics, especially, prior to the invention of writing; and offers a view of the cognitive leaps as mathematics transitioned to the abstract and context-independent phase in three steps - arithmetic, geometry, and algebra.

## 2. Arithmetic

The Ṛg Veda (RV), the most ancient Indian text, displays the knowledge of numbers and arithmetic in around 1500 BCE. It mentions (RV 1.164.10-15, 48) numbers three, five, six, seven, ten, twelve, three-hundred-and-sixty, and seven-hundred-and-twenty. Verse 1.164.11 speaks of a year with 720 days (*ahas*) and nights (*rātris*) joined in pairs, together suggesting a year of 360 days. The two verses RV 10.62.7-8 mention a thousand (*sahasra*) cows and a hundred (*śata*) horses. There are other texts as well, such as *Taittirīya Saṃhitā* sections 7.2.11-20, and the use of poetic meters or *chandas* that offer ample evidence of the knowledge of numbers. Examples of large numbers are not uncommon in these texts (Murthy, 2005; Bavare and Divakaran, 2013).

But these texts neither mention number symbols nor mathematical operations between them. With the absence of any script for writing in the period 1,800 to about 300 BCE, the question is how these numbers were represented and calculations performed.

Interestingly, there is significant evidence of counting, creating number words or number names and some arithmetical operations in that period. The *Śatapatha Brāhmaṇa* (SB; Eggeling, 1897), a middle of the first millennium BCE text, seemingly displays the knowledge of numbers and multiplication: Sections SB 10.4.2.2-17 present fifteen numbers 1, 2, 3, 4, 5, 6, 8, 9, 10, 12, 15, 16, 18, 20 and 24 as the divisors of 720 - the 360 pairs of day and night (represented by bricks) mentioned above. However, in the absence of numerals and operations between them, it would be appropriate to consider that SB actually enumerates ways of distributing or splitting 720 bricks into smaller collections. Each verse has the form: He (*Prajāpati*) made himself four bodies of a hundred and eighty bricks each; … He made himself six bodies of a hundred and twenty bricks each; … He made himself ten bodies of seventy-two bricks each; and finally, …He made himself twenty-four bodies of thirty bricks each (SB 10.4.2.4-17). The repeating phrase "made himself" indicates creating collections of concrete objects, bricks, and counting and adding them; which was the practical way of keeping account of things or objects prior to the invention of written symbols for numbers and abstract procedures for addition and multiplication.

In addition, the *yajuṣamatī* bricks and enclosing-stones also represented day (aha) and night (rātri) and their corresponding 15 *muhūrtas* of time, respectively, making three-hundred-and-sixty each of stones and bricks in an altar (SB 10.4.2.30; 10.5.4.5). These yajuṣamatī bricks referred to half-moons (fortnights), months, and seasons, as well, and *lokamprinā* bricks to the muhūrtas in a year (SB 10.4.3.12). The *nakṣatras* (asterisms) that are time markers in their own way, were also represented by these bricks (SB 10.5.4.5).

Thus, while the yajuṣamatī bricks referred to multiple entities, both stones and bricks referred to muhūrtas. Obviously, this system of representation lacked standardization. The objects had context-specific meanings. Hence, the Śatapatha Brāhmaṇa points to an ancient method of representation and computation. At this stage, arithmetic was more of a "proto-arithmetic," embedded in the context of a social activity of fire sacrifice and tied to concrete entities representing abstract notions of time units, time points, and numbers. This kind of mathematics allowed physical operations of collecting, combining, removing, sorting, separating, distributing, apportioning, and counting of concrete objects in contrast to the mathematics with symbolic operations performed on abstracts numerals. While the former was fixed and embedded in its context with a specific purpose and meaning, the latter is independent of context, permitting freedom for its growth by manipulating written symbols.

Not just the bricks and stones, even sunrise, sunset, and the full- and the new-moon also served as tokens or counters for representing numbers. The annual year-long ritual, *satra*, was performed in which oblations were offered every morning (sunrise) and evening (sunset), every new and full moon (i.e., *parvas*), and at the start of every season and *ayana*, when the sun changes its direction of motion from north to south and vice versa at the two solstices (SB 1.6.3.35-36, Eggeling, 1882; Tilak, 1893). These periodic sacrifices served the same purpose as the tally marks on a stick or a piece of bone in prehistoric times used for time-tracking by following lunar phases (Marshack, 1991). The concept of *ṣaḍaha* was yet another tool developed for time tracking by counting, in which a set of six consecutive days were identified by specific rituals (SB pt. III, p. xxi, Eggeling, 1885; SB pt. V, p. 148, Eggeling, 1894; Dikshit, 1969, p. 22). Ṣaḍaha divides a thirty-day month into five equal parts. Every time the full moon returned a day early on the 29th day instead of 30th, the last day of the ṣaḍaha was dropped to keep in sync with the lunar phase.

Thus, the annual ritual functioned as a concrete calendar by tracking solar and lunar movements. It entailed counting the number of days as one moved through the fortnights, months, seasons and a year, one day at a time, dropping a day here, adding another there, using no more than simple addition and subtraction. Using this method, they reconciled the lunar and solar years, and figured out the length of year between 365 and 366 days.

Names were also used as time markers. Ancient Indians were prolific with generating names: Five names for years in a five-year *yuga* (or cycle); the twelve month names; the dark (*kṛṣṇa*) and the bright (*śukla*) fortnights or *parvas* of a month; the names for 24 parvas in a year; the two parts of each day in the light and dark halves in a month adding up to sixty distinct names of days (ahas) and nights (rātris) in a month; the 30 divisions of a day called muhūrtas, each with its own name in the two fortnights or *pakṣa* in a month; and the subdivision of each muhūrta into fifteen *prati-muhūrtas* with their own specific names (Dikshit, 1969). Also, every fortnight has fifteen *tithis* with their ordinal names: *Prathama* (first), *dvitīya* (second), *tṛtīya* (third), *caturthī* (fourth) and so on until the fifteenth day, called either *pūrṇamāsī* (full moon) or *amāvasyā* (new moon), depending on the bright or the dark half of a month, respectively. Even nakṣatras or star names were used as pointers to a time. The practice of using some of these names is still prevalent in India, especially in astrological and religious contexts.

In short, as exemplified by the cases above, in the absence of script and symbolic representation of numbers and operations, natural events, physical objects, and names representing a thing or time were used to count, and add and subtract to calculate.

## 3. Geometry

Geometry was another field of mathematics which was tied to concrete objects and processes. Unlike Euclidean geometry, it was not founded on definitions and axioms; instead, it emerged as practical geometry from measurement of land, more specifically, in altar construction. In order to achieve specific desired goals, altars of different shapes but of the same size, 7 ½ (square) *puruṣa*, had to be built for fire sacrifice. A puruṣa is a unit of length (= 108 *aṅgulas* or 6 feet 9 inches) but the same name was also used for a square with side one puruṣa long. To add to the challenge of construction, the size of these altars was increased by one square puruṣa at each subsequent performance of the *yajña*, thus, scaling up from 7 ½ to 8 ½ to 9 ½ up to 101 ½ square puruṣa while keeping the same shape. The ancient texts, named *Śulvasūtras* (800 - 300 BCE), are the manuals or practical guides for building these altars and fireplaces using ropes/cords, pegs and bamboo that involve solving geometric problems. Their authors, the *sūtrakāras*, were not the geometricians driven by the need to prove theorems of geometry through deductive reasoning; instead, they were altar makers doing their job following common sense, observation, and intuition. In this process, they discovered certain geometrical truths and principles.

Let us look at some of the geometrical operations sūtrakāras performed while making altars. The first thing was to draw the east-west line on a leveled piece of ground. This line, called *prāchī*, was the line of symmetry for the construction and orientation of altars. The sun's motion and associated shifting direction of shadows of objects on the ground were leveraged to accomplish the goal. A peg was stuck in a measured piece of level ground, and a circle was drawn around it by means of a rope with one end of it tied to the peg. Two points at which the tip of the shadow of the stake touched the circumference at sunrise and sunset were marked, and a line drawn through them was taken as the east-west line (*Kātyāyana Śulba Sūtra*, KS, 1.2; Khadilkar, 1974).

The next step was to draw a perpendicular to this line pointing to the other two cardinal directions north and south. Two pegs were stuck at points A and B equidistant from the point C on the east-west line at which the perpendicular was to be drawn (Fig. 1). Then a rope/cord of desired length was tied to the two pegs and the middle point M on the cord was marked. Holding the marked point M the cord was stretched to the north till it was taut with the point M touching the ground, where point P was marked. Likewise, the cord was stretched in the opposite direction and point Q was marked in the south. A line passing through the points P and Q was perpendicular to the prāchī (KS 1.3). Other variations of this method were also used. Using a similar method, Baudhāyana showed how to draw a square whose corners point in the four intermediate directions (north-east, north-west, etc.) such that the body of the square is parallel to prāchī, hence permitting the construction of altars aligned in the east-west direction (Baudhāyana Śulvasūtras, BS 1.22-28, Thibaut, 1875; Dani, 2010).

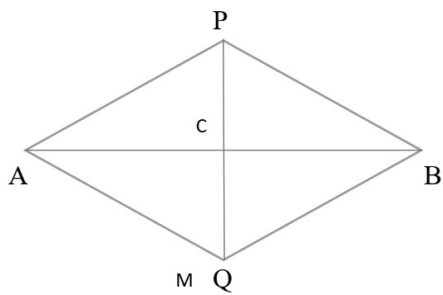

**Figure 1. Drawing a perpendicular line PQ on line AB at point C with a rope**

It is to be noted, that while the sūtrakārs had intuitively grasped the ideas of a straight line and of a circle, they had no concept of angle such as 90-degree angle between two lines perpendicular to each other. The idea of perpendicularity was tied to constructing a line on another line following the method described above. Prāchi was a special case of base-line AB with east-west orientation. Even as late as in the 5th century, the plumb line was the test of verticality (*Āryabhaṭīya* 2.13).

The various constructions in the Sūtras cannot all be presented here. The above example is only to show the method of operation of the sūtrakāras. Once the structure like in Figure 1 was in place and the four points A, B, P and Q on the ground were joined, the structure thus produced was a quadrilateral (*caturasra*) APBQ with its sides AQ=BQ=AP=BP, and AB and PQ its diagonals perpendicular to each other. If the two diagonals were of equal length, it was a square (*samacaturasra*), otherwise, a rhombus (*ubhayataḥprauga*), half of which cut along a diagonal was an Isosceles triangle (*prauga*), such as APB or PBQ. Drawing other shapes, such as trapezia, using ropes and pegs were also possible.

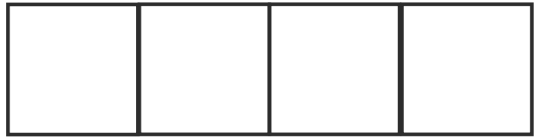

**Figure 2. Forming a rectangle by joining unit squares**

By appending two or more unit-squares (units of any length) side-by-side, a rectangle (*dīrghacaturasra*) could be formed with one side longer than the other (Fig. 2), the area of which could be determined simply by counting the number of unit squares joined together. The 7 square puruṣa altar included seven squares each of area one square puruṣa. Conversely, a rectangle could be divided into unit squares the total number of which would give the area of the rectangle.

Similarly, a square consisted of multiple unit squares. When divided into smaller subunits, it yielded geometrical objects of fractional area. For example, a square with its side 1/2, 1/3ʳᵈ or 1/4ᵗʰ the length of a given measure (unit) was found to have an area 1/4ᵗʰ, 1/9ᵗʰ and 1/16ᵗʰ of the total area produced by the original measure (KS 3.8-10). The opposite of it was also found true: a square with its side twice or thrice or four-times the length of a given unit measured four-fold, nine-fold or sixteen-fold in area, respectively, of the original square with the measure as its side (KS 3.6). Thus, it was established that the area of a square increases or decreases based on the length of its side (measuring cord) (KS 3.12). This conclusion was drawn not by multiplying numbers or fractions, as written above, but by counting the number of unit squares in a given square (KS 3.7).

Following this method of calculating the area of a square, an identity, what is now known as the Pythagoras theorem, was discovered by counting the numbers of unit squares inside the squares drawn on two sides of a rectangle and its diagonal. Notice, this observation was not in reference to the right-angled triangle. Similar observations led to the discovery of Pythagorean triples (3, 4, 5), (5, 12, 13), (7, 24, 25), (8, 15, 17) and (12, 35, 37). And on the basis of these observations, the general conclusion – "in a rectangle, the diagonal produces by itself an area, which its length and breadth produce separately (BS 1.48; KS 2.11)" – was reached by the sūtrakāras, inductively. They had no knowledge of the modern algebraic expression of the theorem $a^2 + b^2 = c^2$, where a and b are the two sides of a right-angle triangle and c its hypotenuse.

Here is another good example of their concrete, hands-on approach, instead of a theoretical one. As they tried measuring the length of the diagonal AB of the square APBQ (assuming AB=PQ in Fig. 1) using AM (=AQ) as the unit length (side of the square), they had to go through several steps of addition and subtraction of pieces of cord (BS 1.61; KS 2.13) to cover the length of AB. First, they stretched the cord AM along the diagonal AB, starting from the point A, the distance MB was left. To cover this gap, they added a third portion of the unit length to the unit AM itself (1 + 1/3), then they further added a fourth portion of 1/3$^{rd}$ of length (1 + 1/3 + 1/3x4), and as this length slightly exceeded the actual length of the diagonal, they subtracted a 34$^{th}$ part of 1/3x4$^{th}$ length, thus covering the length AB, almost. Such small parts may seem indiscernible on a rope, but it depends on the dimension of altars and the length of the rope used. Following an example offered by Dani (2010): For a square of side 35 feet for unit length, the diagonal would be about 50 feet; 1/3$^{rd}$ of the unit (side) length 11.67 feet; 1/12$^{th}$ of the unit about 3 feet (35/3x4 = 2.917 feet); and 1/34$^{th}$ of the last piece would be 0.0858 (=2.917/34) feet or 1.03 inches, which is more than an *aṅgula* (= 3/4$^{th}$ of an inch), making the final correction in 50 feet length manually measurable.

Following these operations, the length of the diagonal found in unit (AM, Fig. 1) lengths was = 1 + 1/3 + 1/3x4 – 1/3x4x34. When expressed in modern decimal notation, this number is equal to 1.414215, which is close to the value 1.4142136 of √2. But this closeness of values should not mislead us into believing that the sūtrakāras performed mathematical calculations by manipulating numbers and fractions or they knew of √2 as an irrational number. They reached the value stepwise through trial-and-error method – by adding and subtracting fractions of a unit length to cover the length of the diagonal. Once the value was established empirically, it became a geometrical fact that the diagonal of a unit square is √2 times the length of its side (unit). This diagonal was called a "doubler," or *dvikaraṇī*, for the square drawn on it as a side was twice in area of the original square (KS 2.12). Note, the sūtrakāras didn't use the modern notation √2 to refer to dvikaraṇī.

Thus, the way arithmetic was tied to bricks and stones and other concrete objects, geometry, based on altar construction, was also concrete and practice bound – using ropes, pegs and bamboo. Extending or reducing the length of a rope by a part of it to draw a circle or a square or other geometrical figures of a certain area is not the same as doing calculations using symbolic numerals. While the former is concrete, the latter is abstract. Since the Brāhmī script appeared in the 3rd century BCE and developed thereafter, it is hard to conceive of the existence of any symbolic manipulation of numbers and fractions prior to that time.

## 4. Algebra

It is with this invention of writing and its development one sees a parallel growth of abstract mathematics in India. *Āryabhaṭa's* (late 5th century) *Āryabhaṭīya* is a testament of such growth in early centuries. The second chapter of the book covers various topics including numbers, fractions, squares and cubes, square roots and cube roots, rule of three, interest calculation, geometry, trigonometry, arithmetic progression series, quadratic and indeterminate equations. While some of these topics are also mentioned in religious, non-mathematical Jain texts such as the *Sthānāṅga Sūtra* from around the turn of the era, absence of any details in them renders it hard to conceive how these calculations were performed, especially when writing was still in its early stages. Even *Arthaśāstra* (3rd century BCE - 3rd century CE) fails to throw light on this matter. It simply mentions metrology – measures of time, length, weight, capacity – and salaries, taxes, tariffs, revenues, rations, etc. Knowledge of fractions is also visible in this text.

This is why the appearance of algebra in Āryabhaṭīya (AB) with no precedence of it in the previous centuries seems puzzling. The lack of direct evidence of it coming from outside the Indian subcontinent also adds to this puzzlement. Hence the questions: Did it develop in India? And how did it develop? Below we shall explore the possibility of its roots in concrete, practical geometry.

As noticed above, geometry was quite prominent in earlier periods. Although mathematicians in the classical period (400-1,600 CE) do not mention altar-making or Śulvasūtras, the knowledge of these Sūtras lasted through later centuries even when altar-making was no longer in practice. Commentaries were written on these Sūtras till as late as the 16th century, and some of them mention the new mathematical works (Ch. 2, Datta, 1932). Sarasvati Amma suggests that the later mathematicians used geometry as an aid to prove the algebraical results (1999, p.220). However, considering the precedence of geometry over algebra in ancient India, the former may also have a role in the inception of the latter.

A path to algebra may be traced through the idea of fractions and equivalence of fractions. Concrete experiences with handling of bricks and ropes and their parts in altar-making led to the knowledge, for example, that one out of two equal parts of a brick (or rope) is the same as two out of four or four out of eight equal parts. This relationship between parts and whole and between numbers made the generalization (in modern symbols) 1/2 = 2/4 = 4/8 = … = a/b possible, where a and b could be any two numbers as long as b is twice as much as a.

This observation found its expression as fractions in symbolic form with the advent of script. Which, then, led to formulation of 'the rule of three' or 'the rule of proportionality' used in solving many problems. Interestingly, Āryabhaṭa mentions 'the rule of three' and manipulation of fractions (ratios) next to each other (see AB 2.26, 27). He uses the rule to solve the gnomon and shadow problems (AB 2.15, 16) in which ratios of perpendicular to base in two similar right triangles are equated. Brahmagupta and Bhāskara II also use the rule for solving geometry problems. Before them, Vedāṅga Jyotiṣa mentions the rule. However, fractional calculations in the extant version of the text were feasible only when it was possible to write and reduce fractions (Dasgupta, 2023). According to Datta and Singh (2004, pp. 188-9), writing of fractions and their reduction came around 200 CE, which is about the time when Brāhmī script and numerals were in use.

Essentially, the rule equates two ratios, such as a/b = x/y, and says that if three quantities (a, b, x) are known, then y, the fourth one, what is called the desired result or *icchāphala*, can be obtained from the relation. The structure of the rule is of the form: that if for quantity a, the result is b, then for quantity x, what will be the result? The symbol y that stands for the unknown quantity or the desired result, may be represented in the form of a simple linear equation with one unknown y = (b/a).x = mx, where m equals b/a. Similar rules of 5, 7, 9 and 11 were apparently derived from the rule of three (Datta and Singh, 2004) to solve more complex problems. These are algorithmic methods, enumerating steps for how to solve practical problems such as calculation of interest on a loan, principal amount, period of loan, or profits or taxes, and so on. Instructions were provided for which quantities in a given problem to be multiplied together and what should be divided by what to find the answer.

Thus, through proportions, though they are now considered part of arithmetic, developed the method for computing unknown values from the known ones in linear equations. They provided general means for solving multiple types of practical problems. These solutions, however, did not use any symbols for unknown quantities. Neither were any operations performed on them. The unknown was only an explicit quantity tied to a concrete something such as the loan amount or the length of a shadow. It was not a variable that required a symbol to represent an unknown quantity as in algebra.

Therefore, while the notions of two characteristics of algebra – one, general solutions, and the other, finding the value of unknown from the known ones – had developed early on with proportions, the latter was still further away from algebra as it had no concept of equation and how to solve them. Unlike in algebra, proportions do not have two sides of an equation that are to be kept balanced. How or when these concepts and techniques appeared in Indian mathematics is not known. But as we shall see below, by the 5[th] century CE mathematicians were aware of them.

An example might help draw the distinction between proportions (arithmetic) and algebra better. Let us consider the problem AB 2.25 in Āryabhaṭīya which is as follows:

A certain amount of loaned money, say P, earned some interest (unknown) in a month. This interest was further loaned at the same rate for T months. The sum A of the original interest and the interest on the original interest is known. What was the original interest earned on P?

Assuming the original interest on P is x per month, and the interest on interest x for the same period of time is y, by applying the rule of three P: x :: x: y, one gets P/x = x/y. Therefore, $y = x^2/P$. Applying the proportionality rule again, i.e., (1: y :: T: ?) if the interest on x for one month is y, then for T months it is equal to $T.y = T(x^2/P) = (T/P).x^2$. If this interest amount were known, it would be possible to find the interest x by first dividing it by T and then taking the square root of yP ($x = \sqrt{yP}$), but instead, the problem gives the sum A (a mixture) of two interests – the original interest x plus the interest on interest $(T/P).x^2$.

So, $$(T/P).x^2 + x = A \qquad \qquad \qquad …(1)$$

This equation adds two dissimilar quantities, or quantities of different dimensions – $x^2$ and x – like adding the area of a square to its side. This is a quadratic equation different from proportions that needed a technique to resolve for the value of the unknown x. Āryabhaṭa provides a formulaic solution without saying how he reached it (Shukla & Sarma, 1976):

"Multiply the interest on the principal plus the interest on that interest by the time and by the principal; (then) add the square of half the principal; (then) take the square root;

(then) subtract half the principal; and (then) divide by the time: the result is the interest on the principal."

In modern notations, the interest $x = 1/T[\sqrt{\{APT + (P/2)^2\}} - P/2]$, which may be computed following above instructions.

Brahmagupta (7[th] century) offers a general solution for quadratic equation (Colebrooke, 1817, p. 346), but it is Bhāskara II (12[th] century) who, in his *Bīja Gaṇita*, first gives the general method for solving quadratic equation (Ibid, pp. 207-9). According to him: both sides of the equation – unknown quantities on one side and known on the other – should be multiplied by an assumed quantity, and some number is to be added on both sides so the unknown side may yield a square root. Equate it with the square root of the known quantity on the other side, and find the value of the unknown quantity. Let us solve equation 1 following the method Bhāskara quotes in BG 5.131 from c. 9[th] century mathematician Śrīdhara (in modern notation):

Step 1: Multiply both sides of equation by a number equal to four times the coefficient $(T/P)$ of $x^2$

$$4.(T/P)^2.x^2 + 4.(T/P).x = 4.(T/P).A$$

Step 2: Then add the square of the original coefficient, 1, of x on both sides

$$4.(T/P)^2.x^2 + 4.(T/P).x + 1 = 4.(T/P).A + 1$$

Simplify the equation by multiplying it with $P^2$ on both sides (not in Śrīdhara):

$$4T^2.x^2 + 4TP.x + P^2 = 4ATP + P^2$$

Therefore,                    $$(2Tx + P)^2 = 4ATP + P^2$$

Step 3: Extract the root [and equate]

$$2Tx + P = +/- \sqrt{\{4ATP + P^2\}}$$

Hence, $x = (1/2T).[\sqrt{\{4ATP + P^2\}} - P] = 1/T.[\sqrt{\{ATP + (P/2)^2\}} - P/2]$          …(2)

Which is exactly the solution Āryabhaṭa provided (the negative value of x is ignored here).

We do not know if he actually followed this method or not. But if he did, then it reflects the knowledge of the algebraic identity $(X + Y)^2 = X^2 + 2XY + Y^2$. Notice the use of this identity in rewriting the expression $4T^2.x^2 + 4TP.x + P^2$ as $(2Tx + P)^2$ above. Although AB 2.25 does not reveal the knowledge of this identity, verse AB 2.23 does. It says: to subtract the sum of squares of two factors (X, Y) from the square of their sum and divide the difference by two to get the product of the two factors, i.e., $X.Y = 1/2[(X + Y)^2 - (X^2 + Y^2)]$, which is just another way of writing the abovementioned identity. With the knowledge of Śulvasūtras still around, it is not hard to see that in Āryabhaṭa's time the identity was known either through geometry (by counting unit squares) or in its algebraic form using symbols. If one draws a square of side $(X + Y)$, one can clearly see the two squares $X^2$ and $Y^2$ in addition to two rectangles of area XY each within the square (Fig. 3). Hence, if Āryabhaṭa followed the procedure above – known as the method of completing the square – he must have used the identity. Indian mathematicians called this method "elimination of the middle term" or *madhymāharaṇa*. The same method was used to solve the problem AB 2.24 where the values of two factors x and y are calculated from the given values of their difference (x - y = a) and their product (x.y = b) that yield a quadratic equation $x^2 - ax = b$. The identity used to solve this equation, $(X - Y)^2 = X^2 - 2XY + Y^2$, too has its geometric representation similar to that in Fig.3.

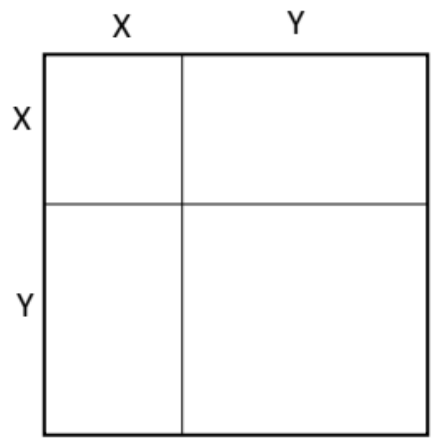

**Figure 3. Geometric representation of identity $(X + Y)^2 = X^2 + 2XY + Y^2$**

The above is a possible scenario of the origin of algebra in India, rooted in arithmetic and geometry. Saying that Indian algebra was rooted in geometry is not the same as one comes across in Babylonian mathematics with geometrical representations of algebraic problems involving not just the measurable lengths and areas but also workers, working days, number of bricks made that were represented by line segments and areas (Høyrup, 2010). In both cases however – Indian as well as Babylonian – the underlying method was the same: by completing the square. While the Babylonians used cut-and-paste geometry method to complete the square, Indians, as one finds in the texts, used arithmetical representation, though some of the concepts and identities were of geometric origin. This origin often remains hidden from sight because of the practical nature of Indian geometry (unlike the Greek geometry with magnitudes only) which dealt with measurable magnitudes or lengths and areas associated with numbers. The modern mindset of viewing algebra with arithmetical representation obscures the geometric roots of it.

As a word of caution, despite the scenario above, the possibility of import of algebra from outside cannot be completely ruled out. It is puzzling to see that while Āryabhaṭa starts his chapter on mathematics with the decimal system (AB 2.2), just ten verses later in AB 2.11-12 switches over to the Babylonian/Greek sexagesimal system of dividing the circle. The 21,600 arc-minutes circumference of a circle is the Babylonian way of dividing a circle into 360 equal parts or degrees each of which is further divided into 60 minutes (360 x 60 = 21,600). Even the models used to describe the planetary motion, though modified in Āryabhaṭīya, are of Greek origin. Hence, assimilation of algebra from outside into Indian mathematics is not impossible. Nevertheless, such absorption could not have been possible without a substantial indigenous knowledge and skill base ready to take it up and apply in the indigenous context. Some of the characteristics of algebra, as noted, had already developed in arithmetic, and the knowledge of certain algebraic identities was present in geometry.

Whichever way it was, in the process of development of algebra, mathematics made its final shift from concrete to abstract. Prior to that, proportions helped solve concrete problems with *icchāphala* as a place holder for the cost of a commodity or the length of a shadow. With them, the notion of the first two of the three characteristics of algebra – generality, unknown quantity, and equation – had already developed. When the idea of equation and how to handle them entered in mathematics isn't quite clear. Bhāskara I (7[th] century) was the first to use the term *samī-karaṇa* or *sama-karaṇa* for equations in the sense of "making equal" in his commentary on verse AB 2.30 (Shukla, 1976, p. 127). Essentially, it says that if a combination of $a_1$ *gulikās* (beads, unknown) and $b_1$ *rūpakas* (coins, known) is *equal in value* to $a_2$ *gulikās* and $b_2$ *rūpakas*, then the value of a bead or the unknown = $[(b_2 - b_1)/(a_1 - a_2)]$ coins. This is a simple linear equation with one unknown (*ek-varṇa-samī-karaṇa*). But in the absence of evidence of use of symbols for equality and unknown, it may be assumed that the meanings were understood from the context in the problem. Bhāskara I, in his *bhāṣya*, used 'o' for unknowns on p. 114 (Ibid., AB 2.25), not '*yāvat-tāvat*' ('as many as' or 'as much as'), while on page 128 he used the word yāvat-tāvat but not the symbol 'o' (AB 2.30). Neither did he set the equation for AB 2.25 nor he showed how to solve it. Instead, he chose to use the formula provided by Āryabhaṭa to find answers for the example problems. It is to be born in mind that knowledge of equality does not imply the knowledge of resolving equations, especially the second and higher order equations. Special techniques are required to solve them. So, how Āryabhaṭa actually set the equations and with what symbols for the known and unknown quantities remains a mystery. All he said was that the knowledge presented in the book was honored at Kusumapura (AB 2.1, Shankar and Shukla, 1976).

## 5. Discussion

The sections above chart the history of mathematics in India over a period of about two millennia during which mathematics passed through several cognitive leaps. One of the major leaps was the movement from concrete to abstract thinking. The mathematician Keith Devlin writes in his book *The Math Instinct* that: "Using our fingers to count indicates that we have a conception of numerosity, but it does not necessarily imply that we have a concept of *number*, which is purely abstract. … [I]f I say 'Think of the number five,' I am using the word 'five' as a noun. As such, it refers to a particular object. What object? *The number five*. The number five is not a concrete object like a chair, but an abstract object. We can't touch it or smell it. But we can think about it. And we can use it." (Italics in the original, Devlin, 2005, p. 206)

When it is said that arithmetic in the Śatapatha Brāhmaṇa was concrete, it does not mean that Indians in the 1st millennia BCE had only a vague sense of number or numerosity. Instead, they had already advanced from the idea of numbers tied to body parts, say a finger or a digit, to the use of objects, such as bricks and stones, and could actually do computation by combining or separating them. Redistribution of 720 bricks into collections of different sizes shows that they had gone beyond counting on fingers to using concrete counters. Not only did a single brick represent a day, it also represented larger quantities, e.g., a fortnight or a month. But this representation was always a one-to-one correspondence, i.e., one brick to a day or a fortnight, a month or a season. Further, these representations were always tied to a pattern in nature. Even when it was a block of 15 muhūrtas, it was essentially a night or a day, never any arbitrary number of muhūrtas, say 8. Likewise, never any number of days other than a fortnight or a month. Hence, the representation was always of a segment of time tied to the motions of the sun or the moon.

In this sense, although arithmetic had freed itself from finger enumeration and hence gained the freedom of counting, adding and subtracting a much larger number of objects than what would have been possible using fingers, it was still tied to natural cycles and their images in the fire altars. Counting of days in the observance of annual ritual and the concept of ṣaḍaha have the same idea: each day was a marker, a notch or a tally to be counted in the cycle of 360/365 days. This is what made arithmetic context-based and concrete - tied to the rituals of altar making and fire sacrifices.

Geometry too had its share of body-parts based length measurements – e.g., aṅgula (width of the middle finger), *hasta* (length of a palm), *aratni*, and puruṣa – but had almost freed itself from it with the use of a concrete object, a rope. A piece of a rope could represent a length, say, a cubit or an aratni, or the height of a man, a puruṣa, in building altars. The use of sesamum seed or *tila* to subdivide an aṅgula into 34 parts is yet another example of measures independent of the human body (Dani, 2010). But geometry was still tied to altar-making. It was not an endeavor discovering geometrical truths separate from rituals.

So, the first cognitive jump came with the separation of numbers and measures from the body parts to concrete objects. It allowed arithmetical calculations and geometrical constructions that would not have been possible otherwise. The use of objects in Śatapatha Brāhmaṇa and the Śulvasūtras for calculation indicates that the separation had occurred before them.

Subdivision of an aṅgula, or two halves of a rectangle, or two bricks representing a full day – one for the daytime (aha), another for the nighttime (rātri) – indicate the next conceptual break in thinking in terms of fractions or parts of a whole. These were unit fractions of which ample examples may be found in the Arthaśāstra. But as remarked above, reduction and calculation of fractions aren't possible without symbolic numerals, though to a limited extent they could be used in practical life. For example, $1/6^{th}$ tax on agricultural production could be paid by separating the sixth measure after every five measures by weight of any unit. The point is that having a rudimentary idea of fractions and its implementation in daily practical life is different from performing mathematical calculations using fractions, the rules for which were given by Āryabhaṭa in the $5^{th}$ century. Some of these rules were possibly known before him.

A critical moment came in the history of Indian mathematics with the rise of the Mauryan empire in the $3^{rd}$ century BCE, when there was push for growing agriculture, irrigation, commerce, taxation, revenue generation, army, weaponry, and construction of roads, forts, palaces, gardens, water reservoirs and so on. Under these conditions, means for record keeping, communication, and calculations for secular projects (free from altar making) were needed. For which the two skills – writing and mathematics – had to be developed as tools for supporting the state administration. It is no surprise to see the emergence of Brāhmī script in this period. A closer look at the tables presented by Datta and Singh (III-XI, Hist. of Hindu Math. Vol 1, 2004) shows both the evolution of and increase in symbolic numerals from the Ashokan period onwards. By the 1st and 2nd century CE there were individual symbols for numbers 1 through 10, multiples of 10, of 100, and even of 1000, which remained in use until 6th or 7th century. Brāhmī neither had a symbol for zero nor was it a system with positional value. Although some composite numbers are presented in table XI, it remains to be determined how mathematical calculations were actually performed using these Brāhmī symbols. Nonetheless, they clearly establish that around the turn of the era Indian mathematics was transitioning from concrete to abstract conception of numbers such that instead of a collection of objects representing a number, say five, a written symbol was sufficient to display the quantity. These numerals were a standard set of arbitrary symbols (though undergoing some changes in shape in the subsequent centuries). And operations could be performed on them that were earlier done manually on a collection of things, e.g., bricks (see above). The use of objects, however, did not disappear immediately after the advent of numerals. The 2nd century CE Buddhist monk Vasumitra mentions the use of counters for calculation (Bronkhorst, 1994).

Assuming calculations with Brāhmī numerals were possible (as there are no examples from that period), the ten topics of mathematics in the Sthānāṅga Sūtra seem feasible, one of which was *kalā-savarṇa* or fractions. With this came the next mental leap of establishing equivalence of ratios, or proportions, of two different things, that allowed, for example, determining the cost of bricks for construction if the price for a given number of them was given, or the calculation of interest on a loan over a duration if the rate of interest was known. But the rule of three (proportions) was not sufficient for solving all types of problems. They could not help when a mixture of quantities of different dimensions, e.g., interest and interest on interest were provided (AB2.25). New techniques were required to solve such problems. Movement from arithmetic to algebra was the next jump.

In addition to the numerals representing known quantities in arithmetic, algebra also needed a mental leap in the use of a symbol to represent an unknown quantity. Although proportions had the notion of an unknown quantity, they did not use symbols for it. The unknown in them was tied explicitly to quantities such as interest or profit. It was never named with a symbol as in algebra. According to Dantzig:

"*The symbol is not a mere formality*; it is the very essence of algebra. Without the symbol the object is a human perception and reflects all the phases under which the human senses grasp it; replaced by a symbol the object becomes a complete abstraction, a mere operand subject to certain indicated operations." (Italics in the original; Dantzig, 2005, p.83)

A symbol thus transcends objects. The way symbol 2 represents any two objects, be they two tables or two persons, the symbol x in algebra represents any number of any object. This x is an operand on which mathematical operations are performed, which was not the practice in proportions wherein operations were performed only on given values. Thus, naming of an unknown and performing operations with it are what distinguish algebra from the numerical problem-solving techniques.

With the above in mind, the symbolic representation leads to transition of literal problems into forms or patterns such as $x^2 + ax - b = 0$, or $x^2 + y^2 = z^2$, or $ax^3 + bx^2 + cx + d = 0$, which show the relationship between known and unknown quantities in a problem. Following Indian symbolism, a quadratic equation such as $10x - 8 = x^2 + 1$ was written as follows (Shukla, 1976, p. lxxvii; Colebrooke, p. 346):

$$\textit{yāva } 0 \textit{ yā } 10 \textit{ rū } \overset{\circ}{8}$$

$$\textit{yāva } 1 \textit{ yā } 0 \textit{ rū } 1$$

Where *yā* is the short form of the unknown quantity *yāvat-tāvat*, *rū* is for *rūpa* or the absolute quantity, and *yāva* is abbreviation of *yāvat-tāvat-varga* meaning square (*varga*) of the unknown quantity *yāvat-tāvat* or *yā* (x). Hence, *yāva* = x². The sign of a circle, as shown above 8 in the equation, was used for negative by Bhāskara I, while the signs of + or x for negative were also used in a few works written in the north-western districts (Shukla, 1976, footnotes, p. lxxii). Although the terms *kṣe* for addition and *gu* for *guṇakāra* or multiplication were used in mathematics, in algebraic equations they were often skipped as in the example above. Likewise, no sign for equality (e.g., =) was used, instead the two sides (or *pakṣa*) of an equation were written one on top of another with their equivalence or balance "understood."

Regarding solving quadratic equations, Shukla suggests that Bhāskara I perhaps knew about the method of completing the square. He interprets the term *yāvakaraṇa* as making the square of a quantity and points to Bhāskara I using the word as a synonym of *varga* or square (-ing). But he also mentions that "nowhere in his commentary has [Bhāskara I] used that term." As observed above, Bhāskara I neither set the equation for AB 2.25 nor did he solve it. Instead, it was Śrīdhara in the 9[th] century who gave the method which was later generalized by Bhāskara II.

Anyway, staying within the scope of this article, it is clear that with the advent of algebra around the middle of the first millennium of this era, Indian mathematics took the final steps in moving away from concrete and context-bound math to context-free and abstract mathematics. This is when the mathematicians started setting equations (based on literal problems) and finding methods for simplifying and solving them for the value of unknown without thinking of the concrete. For instance, in the process of solving the equation for AB 2.25 in Section 4 no thought was engaged with the interest rate. Instead, the procedure involved handling only mathematical objects in symbolic forms. This freedom from concrete permitted further growth of algebra to solving third- and fourth-degree equations as well as indeterminate equations with more than one unknown.

## 6. Conclusion

From a modern perspective, the mathematics discussed above is of elementary level. But this is the kind of math that kept the best mathematicians thinking in ancient India. While studying and cataloguing their accomplishments, it is also important to understand the conditions under which these advances were made because that tells the story of human ingenuity and shows how ideas develop and mathematics progresses. This article is about such progress in ancient India. From growing complexity of social structure and material transactions arose the need for cognitive leaps for solving increasingly complex mathematical problems. It took about a thousand years to transition from concrete arithmetic and geometry to abstract algebra in India. Similar to Babylonian mathematics that progressed from using tokens and bullae to abstract form, Indian mathematics too started its journey from the use of bricks and stones, and ropes and pegs in altar making till it reached algebra. Proportions, prior to algebra, had their origin in the idea of ratio tied to the experience of concrete objects and were used for solving practical problems. Abstract notions of – general solution, unknown quantity, and equation – the three characteristics of algebra, developed as mathematics proceeded from proportions to algebra, seeking methods of solution free from concrete objects and context-dependence. This progress became possible with symbolic representation of numbers and unknown quantities that came with the invention of writing, and freed up mathematics for its future growth. Unlike in Sumerian history, though, so far, there is no evidence showing numbers driving the introduction of symbolic representation in India. Instead, mathematics seems to have leveraged the script as it appeared around 3rd century BCE.

## Acknowledgement

The author would like to express his gratitude to Dr. Vasant K. Nagulapalli and the reviewer for their constructive comments and valuable advice.

## References

[1]    C. Barras. How did ancient humans learn to count? Nature, 594, pp. 22-25, 2021.

[2]    B. Bavare and P. P. Divakaran. Genesis and Early Evolution of Decimal Enumeration: Evidence from Number Names in Ṛgveda. Indian J. Hist. Sci. 48(4), pp. 535-581, 2013.

[3]    J. Bronkhorst. A Note on Zero and the Numerical Place-Value System in Ancient India, Asiatische Studien / Études Asiatiques 48(4), pp. 1039-1042, 1994.

[4]    M. Clagett. Ancient Egyptian Mathematics V3, Am. Phil. Soc., 1999.

[5]    H. T. Colebrooke. Algebra with Arithmetic and Mensuration of Brahmegupta and Bhascara, London, 1817.

[6]    P. Damerow. The Origins of Writing and Arithmetic. In: The Globalization of Knowledge in History, Max-Planck Inst. for the Hist. of Sci., 2012.

[7]    S. G. Dani. Geometry in the Śulvasūtras. In: C.S. Seshadri ed. Studies in the History of Indian Mathematics, pp. 9-37. Hindustan Book Agency: New Delhi, 2010.

[8]    T. Dantzig. Number: The Language of Science. Pi Press: NY, 2005.

[9]    J. Dasgupta. Emergence of Indian Mathematics in the Vedic Period: A Reassessment. 2023. (under review).

[10]   B. Datta. The Science of the Sulba: A Study in Early Hindu Geometry. Univ. Calcutta, 1932.

[11]   B. Datta and A. N. Singh. History of Hindu Mathematics Vol. 1. Bhartiya Kala Prakashan: Delhi, 2004.

[12]   K. Devlin. The Math Instinct: Why You're a Mathematical Genius. Thunder's Mouth Press: New York, 2005.

[13]   S. B. Dikshit. Bhartiya Jyotish Sastra V1. Govt. India Press, 1969.

[14]   J. Eggeling (Tr.). The Śatapatha-Brāhmaṇa. Clarendon Press: Oxford, 1882-97.

[15]   J. Høyrup. Old Babylonian "Algebra", and What it Teaches Us about Possible Kinds of Mathematics. In: Mathematics in Ancient Times. ICM Conference, Kerala School of Mathematics, 2010.

[16]   S. D. Khadilkar (Ed.). *Kātyāyana Śulba Sūtra*. Vaidika Samsodhana Mandal: Poona, 1974.

[17]  A. Marshack. The Roots of Civilization: The Cognitive Beginnings of Man's First Art, Symbol and Notation. Moyer Bell Ltd., 1991.

[18]  S. S. N. Murthy. Number Symbolism in the Vedas. Electronic J. Vedic Studies 12(3), pp. 86-98, 2005.

[19]  T. A. Sarasvati Amma. Geometry in Ancient and Medieval India. Motilal Banarsidass Pbl., Pvt. Ltd.: Delhi, 1999.

[20]  D. Schmandt-Besserat. How Writing Came About. Univ. of Texas Press: Austin, 1996.

[21]  K. S. Shukla (Ed.). Āryabhaṭīya of Āryabhaṭa with the Commentary of Bhaskar I and Someśvara. INSA: New Delhi, 1976.

[22]  K. S. Shukla and K. V. Sarma (Trs. and Eds.). Āryabhaṭīya of Āryabhaṭa. INSA: New Delhi, 1976.

[23]  G. Thibaut. *Śulvasūtra* of Baudhāyana. The Pandit, 1875.

[24]  B. G. Tilak. The Orion. Mrs. Radhabai Atmaram Sagoon Pbl.: 1893.

## Contact Details:

**Jaidev Dasgupta, Ph.D.**

An Independent Researcher who was trained as a scientist at Panjab University, Chandigarh, and Tata Institute of Fundamental Research (T.I.F.R.), Mumbai in India, and worked at Vanderbilt and Harvard universities in USA. His current interests lie in the history of science and mathematics in ancient India.

E-mail: jaidevd101@gmail.com